\documentclass[11pt,leqno]{article} 
\usepackage{graphics}

\let\OLDthebibliography\thebibliography
\renewcommand\thebibliography[1]{
  \OLDthebibliography{#1}
  \setlength{\parskip}{2pt}
  \setlength{\itemsep}{2pt plus 0.3ex}
}

\newtheorem{thm}{Theorem}[section]
\newtheorem{lma}{Lemma}[section]

\newtheorem{prop}{Proposition}

\newcommand{\beqa}{\begin{eqnarray}}
\newcommand{\eeqa}{\end{eqnarray}}

\newcommand{\pf}{\noindent {\bf Proof:} $\s$ }
\newcommand{\epf}{ \hfill$\diamondsuit$ \medskip}

\newcommand{\beq}{\begin{equation}}
\newcommand{\eeq}{\end{equation}}
\newcommand{\lbl}{\label}
\newcommand{\s}{\; \;}

\newcommand{\ep}{\epsilon}

\newcommand{\ra}{\rightarrow}
\newcommand{\al}{\alpha}
\newcommand{\p}{\varphi}

\title{Nonlinear oscillators at resonance with periodic forcing}

\author{
Philip Korman   \\ 
Department of Mathematical Sciences \\ 
University of Cincinnati \\ 
Cincinnati Ohio 45221-0025 \\
and \\
\\
 Yi Li \\
Department of Mathematics  and Statistics\\
John Jay  College of Criminal Justice \\
New York, NY 10019
}

\date{}

\begin{document}

\maketitle
\begin{abstract} 
In this note we unify the results of A.C. Lazer and P.O. Frederickson \cite{FL}, A.C. Lazer \cite{L2}, A.C. Lazer and D.E. Leach  \cite{L1}, J.M. Alonso and R. Ortega \cite{O},  and P. Korman and Y. Li \cite{KL} on periodic oscillations  and unbounded solutions of nonlinear equations with linear part at resonance and periodic forcing. We give conditions for the existence and non-existence of periodic solutions, and obtain a rather detailed description of the dynamics for nonlinear oscillations at resonance, in case periodic solutions do not exist.
\end{abstract}

\begin{flushleft}
\end{flushleft}

\begin{flushleft}
\end{flushleft}

\section{Introduction}
\setcounter{equation}{0}
\setcounter{thm}{0}
\setcounter{lma}{0}

We are interested in the existence of $2 \pi$ periodic solutions to the problem (here $x=x(t)$)
\beq
\lbl{1}
x''+f(x)x'+g(x)+n^2x=e(t) \,.
\eeq
The linear part
\[
x''+n^2x=e(t)
\]
is at resonance, with the null space spanned by $\cos nt$ and $\sin nt$. 
Define $F(x)=\int_0^x f(z) \, dz$. We assume throughout this paper that  $e(t) \in C(R)$ satisfies $e(t+2\pi)=e(t)$ for all $t$, $f(x),g(x) \in C(R)$, $n \geq 1$ is an integer; moreover, we assume that the finite limits at infinity $F(\infty)$, $F(-\infty)$, $g(\infty)$, $g(-\infty)$  exist, and
\beq
\lbl{2}
F(-\infty)<F(x)<F(\infty),  \s\s \mbox{for all $x \in R$} \,,
\eeq
\beq
\lbl{2.1}
g(-\infty)<g(x)<g(\infty),  \s\s \mbox{for all $x \in R$} \,.
\eeq
Define
\[
A_n=\int_0^{2 \pi} e(t) \cos nt \, dt, \s B_n=\int_0^{2 \pi} e(t) \sin nt \, dt \,.
\]

In the case when $f=0$, the equation 
\beq
\lbl{1a}
x''+g(x)+n^2x=e(t) 
\eeq
was considered in the  paper of A.C. Lazer and D.E. Leach \cite{L1} who proved the following classical theorem.

\begin{thm} (\cite{L1} )
The condition 
\[
\sqrt{A_n^2+B_n^2}<2\left(g(\infty)-g(-\infty) \right)
\]
is necessary and sufficient for the existence of $2\pi$ periodic solutions of (\ref{1a}).
\end{thm}

In case $g=0$, for the corresponding equation 
\beq
\lbl{1b}
x''+f(x)x'+n^2x=e(t) 
\eeq 
one has the following result.

\begin{thm}\lbl{thm:2}
The condition 
\[
\sqrt{A_n^2+B_n^2}<2n \left(F(\infty)-F(-\infty) \right)
\]
is necessary  and sufficient  for the existence of $2\pi$ periodic solutions of (\ref{1b}).
\end{thm}

This theorem was proved by A.C. Lazer \cite{L2} for $n=1$ (with an earlier result by P.O. Frederickson and A.C. Lazer \cite{FL}), and for all $n \geq 1$ by P. Korman and Y. Li \cite{KL}, who used a small modification of A.C. Lazer's proof. 
\medskip

Question: is it possible to combine these theorems for the equation (\ref{1})? It turns out that the necessary conditions can be combined, while sufficient conditions cannot be combined.

\begin{prop}\lbl{prop:30}
The condition 
\beq
\lbl{3}
\sqrt{A_n^2+B_n^2}<2n \left(F(\infty)-F(-\infty) \right)+2\left(g(\infty)-g(-\infty) \right)
\eeq
is necessary  for the existence of $2\pi$ periodic solution of (\ref{1}).
\end{prop}

\begin{prop}\lbl{prop:20}
The condition (\ref{3}) is not sufficient for the existence of $2\pi$ periodic solution of (\ref{1}).
\end{prop}

In case   (\ref{1a}) has no $2 \pi$ periodic solutions, all solutions of (\ref{1a}) are unbounded as $t \ra \pm \infty$, as follows by the second Massera's theorem, as was observed first by G. Seifert \cite{s}.  Later J.M. Alonso and R. Ortega \cite{O} gave an elementary approach to this result (with a more refined statement, asserting that solutions tend to infinity in $C^1$ norm). We observe next that the  approach  of \cite{O}  works for the equation (\ref{1}) as well.

\begin{prop}\lbl{prop:1}
Assume in addition to the assumptions above that $f(x)$ is uniformly bounded from below (the assumption (\ref{++*}) below). Then in case 
\beq
\lbl{3.1}
\sqrt{A_n^2+B_n^2} \geq 2n \left(F(\infty)-F(-\infty) \right)+2\left(g(\infty)-g(-\infty) \right),
\eeq
all solutions of (\ref{1}) satisfy  $\lim _{t \ra \pm \infty} \left( x^2(t)+{x'}^2(t) \right)=\infty$.
\end{prop}

Clearly, there are no $2\pi$ periodic solutions in this case, in view of Proposition \ref{prop:30}. Proposition \ref{prop:1} shows that the absence of $2\pi$ periodic solutions turns out to be more decisive in determining the overall dynamics of (\ref{1}) than the existence of $2\pi$ periodic solutions.

\section{The proofs }
\setcounter{equation}{0}
\setcounter{thm}{0}
\setcounter{lma}{0}

The following elementary lemmas are easy to prove.
\begin{lma}\lbl{lma:1}
Consider a function $\cos (nt-\p)$, with an integer $n$ and any real $\p$. Denote $P_c=\{t \in (0,2\pi) \, | \, \cos (nt-\p)>0 \}$ and $N_c=\{t \in (0,2\pi) \,| \, \cos (nt-\p)<0 \}$. Then
\[
\int_{P_c} \cos (nt-\p) \, dt=2, \s\s \int_{N_c} \cos (nt-\p) \, dt=-2 \,.
\]
\end{lma}

\begin{lma}\lbl{lma:2}
Consider a function $\sin (nt-\p)$, with an integer $n$ and any real $\p$. Denote $P_s=\{t \in (0,2\pi) \, | \, \sin (nt-\p)>0 \}$ and $N_s=\{t \in (0,2\pi) \,| \, \sin (nt-\p)<0 \}$. Then
\[
\int_{P_s} \sin (nt-\p) \, dt=2, \s\s \int_{N_s} \sin (nt-\p) \, dt=-2 \,.
\]
\end{lma}

\noindent
{\bf Proof of Proposition \ref{prop:30}}. Given arbitrary numbers $a$ and $b$, one can find a $\delta \in [0,2\pi)$, so that
\[
a \cos nt +b \sin nt=\sqrt{a^2+b^2} \cos (nt-\delta) \,.
\] 
(with $\cos \delta=\frac{a}{\sqrt{a^2+b^2}}$, $\sin \delta=\frac{b}{\sqrt{a^2+b^2}}$.) It follows that 
\beq
\lbl{*}
\frac{A_n}{\sqrt{A_n^2+B_n^2}}\cos nt+\frac{B_n}{\sqrt{A_n^2+B_n^2}}\sin nt =\cos (nt-\delta)\,,
\eeq
for some $\delta \in [0,2\pi)$.
Multiply (\ref{1}) by $\frac{A_n}{\sqrt{A_n^2+B_n^2}} \cos nt$ and integrate, then multiply (\ref{1})  by $\frac{B_n}{\sqrt{A_n^2+B_n^2}} \sin nt$ and integrate, and add the results:
\beq
\lbl{4}
\s\s\s\s \sqrt{A_n^2+B_n^2}= \int_0^{2 \pi} F(x(t))' \cos (nt-\delta) \, dt+\int_0^{2 \pi} g(x(t)) \cos (nt-\delta) \, dt,
\eeq
in view of (\ref{*}).
Using that  $x(t)$ is a $2\pi$ periodic solution,  and Lemma \ref{lma:2}, obtain
\beqa \nonumber
& \int_0^{2 \pi} F(x(t))' \cos (nt-\delta) \, dt=n\int_0^{2 \pi} F(x(t)) \sin (nt-\delta) \, dt \\ \nonumber
&=n\int_{P_s}F(x(t)) \sin (nt-\delta) \, dt + n\int_{N_s}F(x(t)) \sin (nt-\delta) \, dt \\
& <2n \left(F(\infty)-F(-\infty) \right) \,. \nonumber
\eeqa
Similarly, using Lemma \ref{lma:1} 
\beqa \nonumber
& \int_0^{2 \pi} g(x(t)) \cos (nt-\delta) \, dt<g(\infty) \int_{P_c} \cos (nt-\delta) \, dt+g(-\infty) \int_{N_c} \cos (nt-\delta) \, dt \\  \nonumber
& =2\left(g(\infty)-g(-\infty) \right) \,.  \nonumber
\eeqa 
The condition (\ref{3}) follows.
\epf

\noindent
{\bf Proof of Proposition \ref{prop:20}}.
Consider the equation
\beq
\lbl{n1}
x''+f(x)x'+g(x)+n^2x=E \cos nt,
\eeq
with a parameter $E$. Calculate $A_n=\int_0^{2 \pi} E \cos^2 nt 
\, dt=E \pi$, $B_n=\int_0^{2 \pi} E \cos nt \sin nt \, dt=0$, and $\sqrt{A_n^2+B_n^2}=E \pi$.
Choose $E$ so that 
\beq
\lbl{48}
\s\s E \pi=\sqrt{A_n^2+B_n^2}=2n \left(F(\infty)-F(-\infty) \right)+2\left(g(\infty)-g(-\infty) \right)-\epsilon,
\eeq
with $\epsilon>0$ small, to be specified. The condition (\ref{3}) holds for the equation (\ref{n1}). If this condition were sufficient, we would have a $2\pi$ periodic solution of (\ref{n1}), and  hence  
\beq
\lbl{n2}
x''+F(x)'+n^2x=E \cos nt -g(x) \equiv \bar e(t)\,.
\eeq
Calculate the coefficients $A_n,B_n$ for (\ref{n2}): 
\beqa 
\lbl{n3}
& \bar A_n=\int _0^{2 \pi} \bar e(t)\cos nt 
\, dt=E \pi -\int _0^{2 \pi} g(x) \cos nt \, dt,\\ \nonumber
& \bar B_n=\int _0^{2 \pi} \bar e(t)\sin nt 
\, dt= -\int _0^{2 \pi} g(x) \sin nt \, dt. \nonumber
\eeqa
Since (\ref{n2}) is solvable, by Theorem \ref{thm:2} we have
\beq
\lbl{n4}
\sqrt{\bar A_n^2+\bar B_n^2}<2n \left(F(\infty)-F(-\infty) \right).
\eeq
For any $\ep >0$ we can choose an index $n_0$, so that for $n \geq n_0$
\beq
\lbl{n5}
|\int _0^{2 \pi} g(x) \sin nt \, dt|<\ep,
\eeq
as follows by well known results on oscillatory integrals, see e.g., O. Costin et al \cite{Co}.
Using (\ref{n3}), (\ref{n5}), followed by (\ref{48}), obtain
\beqa \nonumber
& \sqrt{\bar A_n^2+\bar B_n^2}>|\bar A_n|>E \pi -\ep \\ \nonumber
& =2n \left(F(\infty)-F(-\infty) \right)+2\left(g(\infty)-g(-\infty) \right)-2\epsilon \\\nonumber
& > 2n \left(F(\infty)-F(-\infty)\right), \nonumber
\eeqa
contradicting (\ref{n4}), provided we fix $\ep<g(\infty)-g(-\infty)$. 
\epf

We shall prove a generalization of Proposition \ref{prop:1} after several preliminary results. By an obvious modification of its proof, one obtains the following  generalization of Proposition \ref{prop:30}.
\begin{prop}\lbl{prop:2}
Assume that the functions $F(x) $ and $g(x) $ have finite infimums and supremums on $(-\infty,\infty)$. Then the condition
\[
\sqrt{A_n^2+B_n^2}<2n \left(\sup F-\inf F \right)+2\left(\sup g-\inf g  \right)
\]
is necessary  for the existence of $2\pi$ periodic solution of (\ref{1}).
\end{prop}

We shall use the following result that is included in  J.M. Alonso and R. Ortega \cite{O}.

\begin{prop} (\cite{O}) \lbl{prop:o}
Let $G(\zeta,\eta): R^2 \ra R^2$ be a continuous vector function, and let $V(\zeta,\eta): R^2 \ra R$ be a continuous function. With $\xi \equiv (\zeta,\eta)$ consider a sequence
\[
\xi _{n+1}=G \left(\xi _{n} \right), \s\s n \geq 0 \,,
\]
beginning with an arbitrary vector $\xi _0$. Assume that 
\beq
\lbl{+}
V \left( G(\xi)\right)>V(\xi) \,, \s\s \forall  \xi \in R^2 \,.
\eeq
Then $\lim _{n \ra \infty} ||\xi _{n}||=\lim _{n \ra \infty} \sqrt{\zeta _n^2+\eta _n^2}=\infty$.
\end{prop}

\pf
If the sequence $\{ ||\xi _{n}|| \}$ fails to tend to infinity, then $\{ \xi _{n} \}$ has a finite accumulation point $\xi ^* \in R^2$. Let $\{ \xi _{n_k} \}$ be a subsequence tending to $\xi ^*$, with $n_1<n_2< \cdots$. Since $V(\xi _{n+1})=V(G \left(\xi _{n}) \right)>V(\xi _{n})$ by (\ref{+}), the sequence $V(\xi _n)$ is increasing. Then one has
\beqa \nonumber
& V(G(\xi ^*))= \lim _{k \ra \infty} V(G(\xi  _{n_k}))= \lim _{k \ra \infty} V(\xi  _{n_k+1})\\ \nonumber
&  \leq \lim _{k \ra \infty} V(\xi  _{n_{k+1}})=V(\xi ^*) \,,\nonumber
\eeqa
contradicting (\ref{+}). (Observing that $n_k+1 \leq n_{k+1}$.)
\epf

The next lemma says that for solution of (\ref{1}), $x^2(t)+{x'}^2(t)$ cannot increase too much over an interval of length $2\pi$.

\begin{lma}\lbl{lma:*}
Assume assume that $e(t) \in C(R)$ is $2\pi$ periodic, the condition (\ref{2.1}) holds, and moreover assume that
\beq
\lbl{++*}
f(x) \geq \al, \s\s \mbox{for some $\al \in R$, and all $x \in R$}.
\eeq
Then for any initial data $(x(0),x'(0))$, with $c_0=x^2(0)+n^2{x'}^2(0)$, there is a number $c=c(c_0)$ so that the corresponding  solution of (\ref{1}) satisfies
\[
x^2(t)+n^2{x'}^2(t) \leq c, \s \mbox{for all $t \in [0,2\pi]$}.
\]
\end{lma}

\pf
Consider the ``energy" $E(t)=\frac{1}{2}{x'}^2(t)+\frac{1}{2}n^2x^2(t)$. Since
\[
E'(t)=-f(x){x'}^2-g(x)x'+e(t)x'.
\]
By our conditions obtain $E'(t) \leq c_1E(t)+c_2$, and the proof follows.
\epf

We now prove the following generalization of Proposition \ref{prop:1}.

\begin{prop}\lbl{prop:4}
Assume that the functions $F(x) $ and $g(x) $ have finite infimums and supremums on $(-\infty,\infty)$ Assume also that (\ref{++*}) holds. 
In case 
\[
\sqrt{A_n^2+B_n^2} \geq 2n \left(\sup F-\inf F \right)+2\left(\sup g-\inf g  \right) ,
\]
all solutions of (\ref{1}) satisfy  $\lim _{t \ra \pm \infty} \left( x^2(t)+{x'}^2(t) \right)=\infty$.
\end{prop}

\pf
Following J.M. Alonso and R. Ortega \cite{O}, we shall use Proposition \ref{prop:o}. Given $\xi=(\zeta,\eta) \in R^2$, denote by $x(t,\xi)$ the solution of (\ref{1}) satisfying $x(0)=\zeta$, $x'(0)=\eta$. Define the map $R^2 \ra R^2$ by $G(\xi)=\left(x(2\pi,\xi), x'(2\pi,\xi) \right)$, and we shall show that the sequence of iterates
\[
\xi _{n+1}=G(\xi _n) \,, \s n=0,1,2, \ldots
\]
is unbounded for any $\xi _0$. With $\delta$ as defined by (\ref{*}), define the function
\[
V(\xi)=\eta \cos \delta-n \zeta \sin \delta+F(\zeta) \cos \delta \,.
\]
Multiply (\ref{1}) by $\frac{A_n}{\sqrt{A_n^2+B_n^2}} \cos nt$ and integrate, then multiply (\ref{1})  by $\frac{B_n}{\sqrt{A_n^2+B_n^2}} \sin nt$ and integrate, and add the results. In view of  (\ref{*}) obtain as above
\beqa
\lbl{15}
& \int _0^{2 \pi} (x''+n^2x) \cos(nt-\delta) \, dt+\int _0^{2 \pi} F(x(t))'\cos(nt-\delta) \, dt \\
& +\int _0^{2 \pi} g(x(t))\cos(nt-\delta) \, dt =\sqrt{A_n^2+B_n^2} \,. \nonumber
\eeqa
Integrating by parts, we express the first term on the left as
\beqa \nonumber
& \left[x'(t) \cos(nt-\delta) \right] |_{_0}^{2 \pi}+n \left[x(t) \sin(nt-\delta) \right] |_{_0}^{2 \pi} \\ \nonumber
& = \left[x'(2 \pi)-\eta \right] \cos \delta-n \left[ x(2\pi)-x(\zeta)\right] \sin \delta\,,
\eeqa
and the second term as
\[
\left[ F(x(2\pi)) -F(x(0)) \right] \cos \delta +n \int _0^{2 \pi} F(x) \sin(nt-\delta) \, dt\,.
\]
We combine the non-integral terms in (\ref{15}) as $V(G(\xi))-V(\xi)$.
Then (\ref{15}) gives
\beqa \nonumber
& V(G(\xi))-V(\xi)=\sqrt{A_n^2+B_n^2} -n \int _0^{2 \pi} F(x(t)) \sin(nt-\delta)\, dt  \\ \nonumber
&   -\int _0^{2 \pi} g(x(t))\cos(nt-\delta) \, dt   \\  \nonumber
& >   \sqrt{A_n^2+B_n^2}  -2n \left(\sup F-\inf F \right)-2\left(\sup g-\inf g  \right)  \geq 0 \,.   \nonumber
\eeqa
(The first inequality is strict because the functions  $g$ and $F$ are non-constant by (\ref{2}) and (\ref{2.1}).) Hence, the condition (\ref{+}) holds, and Proposition \ref{prop:o} applies, proving  the unboundness of the sequence $\left(x(2n\pi,\xi), x'(2n\pi,\xi) \right)$. If there was a sequence $\{ t_k \} \ra \infty$ with bounded ${x'}^2(t_k)+x^2(t_k)$, we would obtain a contradiction with Lemma \ref{lma:*}.
\epf

\end{document}